\documentclass[10pt]{article}

\usepackage{amsmath,amssymb,dsfont,slashed,color,hyperref,amsthm}
\usepackage{authblk}

\usepackage[toc,page]{appendix}
\usepackage{graphicx}
\usepackage[shortlabels]{enumitem}

\definecolor{armygreen}{rgb}{0.29, 0.33, 0.13}

\newcommand{\R}{\mathbb{R}}
\newcommand{\C}{\mathbb{C}}

\newcommand{\id}{\mathbb{I}}

\newcommand{\im}{{\rm i }}

\newcommand\be{\begin{eqnarray}}
\newcommand\ee{\end{eqnarray}}

\newtheorem*{theorem}{Theorem}
\newtheorem*{corollary}{Corollary}
\newtheorem{definition}{Definition}
\newtheorem{proposition}{Proposition}

\newcommand{\calA}{\mathcal{A}}

\newcommand{\cecs}{Centro de Estudios Cient\'{\i}ficos (CECs), Arturo Prat 514, Valdivia, Chile}
\newcommand{\uss}{Universidad San Sebasti\'an, General Lagos 1163, Valdivia, Chile}
\newcommand{\nott}{School of Mathematical Sciences, University of Nottingham, Nottingham, NG7 2RD, UK}

\begin{document}

\title{Variations on a theme of MacDowell-Mansouri}

\author[1,2]{P. D. Alvarez \thanks{E-mail: \href{mailto:pedro.alvarez@uss.cl}{\nolinkurl{pedro.alvarez@uss.cl}}}}

\author[3]{K. Krasnov \thanks{E-mail: \href{mailto:kirill.krasnov@nottingham.ac.uk}{\nolinkurl{kirill.krasnov@nottingham.ac.uk}}}}
\affil[1]{\uss}
\affil[2]{\cecs}
\affil[3]{\nott}

\maketitle

\begin{abstract}
Inspired by the MacDowell–Mansouri formulation of four-dimensional General Relativity, we study a class of four-dimensional gauge-theoretic functionals obtained from the Pontryagin density of a G-connection by inserting, under the trace, a matrix that breaks the gauge group G to a subgroup H. Concretely, we study the model with the pair $(G,H)$ given by $({\rm SU}(3), {\rm U}(2))$.  We show that the critical points of the resulting functional are constant scalar curvature almost-K\"ahler 4-manifolds. On compact 4-manifolds, a stronger conclusion holds under the additional assumption that the scalar curvature is non-negative and the first Chern class is such that an Einstein metric can exist. In this case results in the literature imply that the critical points are K\"ahler-Einstein 4-manifolds. 
 \end{abstract}

\tableofcontents

\section{Introduction}\label{intro}

In \cite{MacDowell:1977jt} MacDowell and Mansouri have described $\Lambda\not=0$ General Relativity (GR) as a gauge theory of a ${\rm SO}(1,4)$ (or ${\rm SO}(2,3)$, depending on the sign of the cosmological constant $\Lambda$) connection. Such a connection can be interpreted as a Cartan connection, see \cite{Wise:2006sm}, combining the spin ${\rm SO}(1,3)$ connection $w$ together with the co-frame field $e$, schematically 
\be
\label{cartan-conn-intr} 
\mathcal{A} = w + \frac{1}{l} e, 
\ee
where $|\Lambda|=1/l^2$. If we decompose
\[
\mathfrak{g} = \mathfrak{h} \oplus \mathfrak{p}, \qquad \mathfrak{g}=\mathfrak{so}(1,4), \quad 
\mathfrak{h}= \mathfrak{so}(1,3), \quad \mathfrak{p}=\R^{1,3},
\]
then the $\mathfrak{h}$ and $\mathfrak{p}$ parts of the curvature 
\[
\mathcal{F} = d\mathcal{A} + \mathcal{A}\wedge \mathcal{A}
\]
encode the curvature and torsion of $\omega$ respectively
\[
\mathcal{F}\Big|_{\mathfrak{h}} = dw + w\wedge w - \frac{1}{l^2} e\wedge e, \qquad
\mathcal{F}\Big|_{\mathfrak{p}} = d_w e. 
\]
When the torsion vanishes, $w$ is the unique metric torsion-free connection for the metric determined by the co-frame $e$, and the vanishing of the $\mathfrak{h}$ part of the curvature $\mathcal{F}$ is the statement that the Riemann curvature is constant. 

Even more impressively, one can rewrite the Palatini action for GR as a functional of the curvature $\mathcal{F}$ of the Cartan connection $\mathcal{A}$. This is most elegantly done by passing to the spinor representation ${\rm Spin}(1,4)$. Let $\gamma^{1,\ldots,5}$ be the $4\times 4$ $\gamma$-matrices realising the Clifford algebra ${\rm Cl}(1,4)$. Then $\mathfrak{so}(1,4)$ is generated by products of two distinct $\gamma$-matrices. The curvature $\mathcal{F}$ is encoded as a $4\times 4$ matrix of 2-forms. The MacDowell-Mansouri (MDM) action is
\[
S[\mathcal{A}] = \int {\rm Tr}(\gamma^5 \mathcal{F}\wedge \mathcal{F}).
\]
If not for the insertion of $\gamma^5$, this would be the Pontragin density for the connection $\mathcal{A}$. The insertion of $\gamma^5$ breaks the invariance to the subgroup ${\rm Spin}(1,3)$ that commutes with $\gamma^5$, and makes the above action non-topological. Expanding the action in terms of the $\mathfrak{h}$ and $\mathfrak{p}=\mathfrak{g}/\mathfrak{h}$ parts of the connection one recovers the Palatini action for GR (plus a topological term). The same idea can be used to rewrite the action of $N=1, D=4$ anti-De Sitter supergravity, see \cite{MacDowell:1977jt}. Generalisations of this construction to higher-dimensional gravities and supergravities have also been considered, see e.g. \cite{Chamseddine:1990gk} for topological theories (i.e. theories without explicit breaking of the $G$-invariance of the action), and e.g. \cite{Castellani:2017vbi} for a construction of the $D=12$ supergravity action in direct parallel with the MDM construction. 

The purpose of this paper is to initiate exploration of analogues of the MDM construction for more general Cartan geometries. Following \cite{Sharpe}, we recall the following definition
\begin{definition}
A Cartan geometry on $M$ modelled on $G/H$ is a pair $(P,\mathcal{A})$ consisting of a principal $H$-bundle $\pi: P\to M$ over $M$, and a $\mathfrak{g}$-valued 1-form $\mathcal{A}\in \Omega^1(P,\mathfrak{g})$ on $P$ called the Cartan connection, such that
\begin{enumerate}
\item {\bf Soldering condition:} For every $p\in P$ the map $\mathcal{A}_p: T_p P \to \mathfrak{g}$ is a linear isomorphism.
\item {\bf Equivaraince:} For all $h\in H$ we have $(R_h)^* A = {\rm Ad}(h^{-1}) A$.
\item {\bf Reproduction of fundamental fields:} For every $\xi\in \mathfrak{h}$, the fundamental vector field $\xi^\sharp\in TP$ generating the right $H$ action on $P$ satisfies $\mathcal{A}(\xi^\sharp) = \xi$.
\end{enumerate}
\end{definition}
The last two properties are generic for a connection on the principal $H$-bundle. What makes Cartan connections unusual is that $\mathcal{A}$ is $\mathfrak{g}$-valued, not $\mathfrak{h}$-valued, and that one imposes the condition that the tangent space to the $H$-bundle $P$ is isomorphic (at every point) to the Lie algebra $\mathfrak{g}$. In particular this means that the tangent space $TM$ is isomorphic to $\mathfrak{g}/\mathfrak{h}$. We will also assume the reductive property, which is that under the action of $H$ on $\mathfrak{g}$ the latter decomposes into the sum of two $H$-modules
\[
\mathfrak{g} = \mathfrak{h} \oplus \mathfrak{p},
\]
or, equivalently ${\rm Ad}(H) \mathfrak{p}\subset \mathfrak{p}$. In a gauge, the soldering condition means that we can identify the $\mathfrak{h}$ part of $\mathcal{A}$ with an $H$-connection on $M$, and $\mathfrak{p}$-valued part of $\mathcal{A}$ with the soldering form on $M$, exactly as in (\ref{cartan-conn-intr}). 

Let us assume that the dimension of $M$, and thus of $\mathfrak{p}$ 
\[
{\rm dim}(\mathfrak{p}) = 2k
\]
is even. We will be interested in the functionals of $\mathcal{A}$ (and thus of its $\mathfrak{h}$ and $\mathfrak{p}$ parts that we denote by $w, e$ respectively) that are constructed from the curvature 2-form $\mathcal{F}=d\mathcal{A} + \mathcal{A}\wedge \mathcal{A}$ and are of the form
\be\label{functional-intr}
S[\mathcal{A}] = S[w, e] = \int_M {\rm Tr}( \Gamma_1 \mathcal{F} \Gamma_2 \mathcal{F} \ldots \Gamma_k \mathcal{F}).
\ee
Here the trace is taken in some representation $V$ of $\mathfrak{g}$, and $\Gamma_{1,\ldots,k} \in {\rm End}(V)$ are some suitable matrices, which are to be chosen in such a way that the functional is $H$-invariant. If not for the insertion of $\Gamma$'s under the trace, the functional would be topological (the integrand would be a total derivative, or, in cases of some $\mathfrak{g}$ and $n$ combination, zero). However, after the insertion of $\Gamma$'s this is no longer the case and one obtains interesting $H$-invariant functionals of the $H$-connection $\omega$ and the soldering form $e$. When ${\rm dim}(M)$ is odd, one could consider the similar construction but this time with the Chern-Simons functionals of $\mathcal{A}$, inserting a suitable $H$-invariant matrix $\Gamma$'s under the trace. Another comment is that one obtains interesting functionals already by inserting a single $\Gamma$ under the trace, but it is important to keep things general to understand the ambiguities involved in the construction. 

On could consider an even more general class of functionals, by following the Chern-Weil construction of characteristic classes and considering all possible $G$-invariant polynomials in the Lie algebra $\mathfrak{g}$. This gives the building blocks for the construction of the action, as a family of topological action functionals is then obtained by taking all wedge products of characteristic classes that give the top degree form. One can then start deforming the forms giving the characteristic classes by changing the $G$-invariant polynomials on $\mathfrak{g}$ into only $H$-invariant ones. In (\ref{functional-intr}) we consider only the simplest example of such a construction, where a single $G$-invariant polynomial is used, given by the trace of an appropriate power of the curvature 2-form. And the passage from $G$-invariant polynomial to an $H$-invariant one is achieved by the insertion of a number of copies of $\Gamma$. It is clear that this can be generalised, but for the purposes of this paper it will be sufficient to work with functionals of the type (\ref{functional-intr}).

Our interest in the functionals (\ref{functional-intr}) is motivated by the fact that they provide a natural principle for building $H$-invariant functionals for $H$-structures encoded in the soldering forms $e$, together with an $H$-connection $\omega$. We recall that a coframe on an $n$-manifold $M$ is a linear isomorphism $u: T_x M \to \R^n$. The set $F(M)$ of soldering forms at the points of $M$ is naturally a principal (right) ${\rm GL}(n,\R)$ bundle over $M$, with projection $\pi: F(M)\to M$. There is a tautological 1-form $e\in \Omega^1(F(M),\R^n)$ on $F(M)$ defined by the formula
\be\label{soldering-form}
e(v) = u(\pi_*(v)), 
\ee
where $v\in T F(M)$ is an arbitrary vector field on the total space of the bundle and $\pi_*$ is the push-forward under the projection map. The soldering form is ${\rm GL}(n,\R)$ equivariant $\forall g\in {\rm GL}(n,\R)$ we have $(R_g)^* e = g^{-1} e$. We can now define
\begin{definition}
An $H$-structure on $M$ is a principal $H$-subbundle $P\subset F(M)$ of the principal ${\rm GL}(n,\R)$ bundle $F(M)$ of frames of $TM$. Equivalently, it is a reduction of the structure group of $TM$ from ${\rm GL}(n,\R)$ to $H$. It comes equipped with an $H$-equivariant $(R_h)^* e = h^{-1} e$ soldering form $e \in \Omega^1(P, \R^n)$ which is horizontal $e(\xi^\sharp)=0, \forall \xi\in \mathfrak{h}$ and such that for each $p\in P$ the map $e_p: T_p P\to \mathfrak{p}$ descends to a fiberwise linear isomorphism $T_x M \sim P_x \times_H \mathfrak{p}, x=\pi(p)$. 
\end{definition}
Comparing this with the definition of the Cartan connection, we see that the $\mathfrak{p}$ part $e$ of a Cartan connection $\mathcal{A}$ is an $H$-structure on $M$. It is also clear that the $\mathfrak{h}$ part $w$ of $\mathcal{A}$ is a principal $H$ connection on $M$. They are combined together into the Cartan connection $\mathcal{A}$. To put it differently, a Cartan geometry modelled on $G/H$ consists of an $H$-structure on $M$ together with a principal $H$-connection on $M$, assembled into a $\mathfrak{g}$-valued 1-form $\mathcal{A}$. The important subtlety here is that not every $H$-structure gives rise to such a Cartan connection, it is very important that $\mathfrak{p}\sim \R^n$ combines together with $\mathfrak{h}$ into the Lie algebra $\mathfrak{g}$. We can call $H$-structures that give rise to Cartan connections $H$-structures of type $(G,H)$. 

The functional of the type (\ref{functional-intr}) is then an $H$-invariant functional of an $H$-structure of type $(G,H)$ encoded by the soldering form $e$, and an $H$-connection $w$. It has the property that the model space $M=G/H$ is automatically the critical point. Indeed, the Euler-Lagrange equations that result from (\ref{functional-intr}) are
\[
d_\mathcal{A} (\sum_{\text{cyclic} \,\,\,\sigma} \Gamma_{\sigma(1)} \mathcal{F} \Gamma_{\sigma(2)}  \mathcal{F} \ldots \Gamma_{\sigma(k)}  ) =0.
\]
These are clearly satisfied by $\mathcal{F}=0$, which describes the model geometry $M=G/H$. So, the homogeneous space $G/H$ is a critical point of any of the functionals (\ref{functional-intr}), but there are other critical points. 

In this paper we will study only a single example of the above construction, leaving other possible interesting examples to further publications. Our example is $({\rm SU}(3),{\rm U}(2))$ in four dimensions. Thus, we shall study functionals of an almost-Hermitian structure in dimension four. 

It is well-known \cite{HH} that an almost-complex structure on a 4-dimensional manifold $M$ exists if an only if there is a class $c_1\in H^2(M,\mathbb{Z})$ such that $c_1= w_2(TM) \, ({\rm mod}\,\, 2)$ (where $w_2$ is the second  Stiefel-Whitney class) and 
\be
c_1^2 = 2\chi(M) + 3\tau(M),
\ee
where $\chi(M),\tau(M)$ are the Euler characteristic and signature of $M$. For example, there is no almost complex structure on the 4-sphere $S^4$ because for it $c_1=\tau=0$, but $\chi=2$. 

If an almost-complex structure on $M$ exists, let $(g,J)$ be an almost-Hermitian structure, which is a pair consisting of a Riemannian metric on $M$, as well as a compatible (orthogonal) almost-complex structure $J:TM\to TM, J^2=-\mathbb{I}$. Let $\alpha,\beta\in \Lambda^{1,0}$ be an orthonormal basis of $(1,0)$ forms, so that the metric takes the form (\ref{metric}). We can encode these into a $\C^2$-valued 1-form 
\[
\Psi = \left( \begin{array}{cc} \alpha \\ \beta \end{array}\right).
\]
Thus, $\Psi$ is a section of the ${\rm U}(2)$ frame bundle $P\to M$ over $M$. Let $w$ be a connection on $P$. The Cartan connection is a sum 
\[ \mathcal{A} = w + \Psi. \]
Given that $\mathfrak{su}(3)=\mathfrak{u}(2)\oplus \C^2$, the Cartan connection is $\mathfrak{su}(3)$-valued. 

We note that we can further split a ${\rm U}(2)$ connection $w$ on $P\to M$ into its traceless and trace parts
\be\label{A-a}
w= A+ a.
\ee
Then the trace part $a$ induces a ${\rm U}(1)$ connection on the determinant line bundle ${\rm det}(TM,J)= K^{-1}$ (anti-canonical bundle). In the used anti-Hermitian conventions where $a$ is purely imaginary, its curvature represents the first Chern class
\be
\left[ \frac{\im da}{2\pi} \right] = c_1(TM,J).
\ee

The MDM construction (\ref{functional-intr}) gives a ${\rm U}(2)$-invariant functional for the resulting ${\rm U}(2)$ structure. We summarise our results in a series of propositions.
\begin{proposition} The most general MDM functional for a pair $(w,\Psi)$ of an ${\rm U}(2)$ connection $w$ and a ${\rm U}(2)$ frame $\Psi$ is given by
\be\label{lagr-intr}
S[w,\Psi]= \int_M {\rm Tr}\Big( \lambda  \Psi^\dagger F\Psi +  \mu da (\Psi^\dagger \Psi) + \nu  (\Psi^\dagger \Psi)^2\Big),
\ee
with real parameters $\lambda,\mu,\nu$ satisfying
\be\label{relation-intr}
3\lambda - \mu + 4\nu=0.
\ee
Here $A,a$ are the ${\rm SU}(2)$ and ${\rm U}(1)$ parts of $w$ as in (\ref{A-a}), and $F=dA+A\wedge A$ is the  curvature of the ${\rm SU}(2)$ connection $A$.
\end{proposition}

\begin{proposition} The Euler-Lagrange equation arising from (\ref{lagr-intr}) by varying with respect to the ${\rm SU}(2)$ connection $A$ is algebraic. Its solution is $A={\bf w}$, where ${\bf w}$ is the anti-self dual part of the Levi-Civita connection. After $A={\bf w}$ is substituted back into the Lagrangian (\ref{lagr-intr}), one obtains a functional of ${\rm U}(2)$ structure encoded by $\Psi$, and the remaining ${\rm U}(1)$ part $a$ of $w$. This functional can be rewritten as that of a pair $(g,\omega)$ of a metric $g$ given in terms of $\alpha,\beta$ by (\ref{metric}) and an almost-symplectic 2-form $\omega$ given by (\ref{omega}), and is given by
\be\label{lagr-gw-intr}
S[g,\omega,a] = -  \lambda \int_M dV ( s + 2\tilde{\mu}-6) -2\tilde{\mu} \omega \wedge \im da .
\ee
Here $dV$ is the volume form, $s$ is the scalar curvature, we have defined $\tilde{\mu}=\mu/\lambda$ and used the relation (\ref{relation-intr}). Viewed as a functional of $(g,\omega)$, these are subject to the compatibility condition 
\be\label{compat-intr}
\omega_{\mu\alpha} g^{\alpha\beta} \omega_{\beta\nu} = - g_{\mu\nu}.
\ee
\end{proposition}

\begin{proposition} The critical points of (\ref{lagr-gw-intr}) are triples $(g,\omega,a)$ satisfying
\be\label{feqs-intr}
Ric^{(2,0)}=0, \qquad da^{(2,0)}=0, \qquad d\omega=0,
\ee
as well as
\be\label{main-feq-intr}
\rho = \frac{1}{2}\omega (s -6)  +\tilde{\mu} (\omega - \im da).
\ee
The first equation in (\ref{feqs-intr}) can alternatively be stated as saying that the Ricci tensor is $J$-invariant. This implies that $\rho(X,Y)= Ric(X,JY)$ is a $(1,1)$ form, and we wrote the last equation in terms of this 2-form. Note that because of the second equation in (\ref{feqs-intr}) all terms in the equation (\ref{main-feq-intr}) are in $\Lambda^{1,1}$. 
\end{proposition}

Using the known result \cite{Draghici} that the Ricci form of an almost-K\"ahler four-dimensional metric with a $J$-invariant Ricci tensor is closed, together with the other equations in (\ref{feqs-intr}), (\ref{main-feq-intr}) we obtain our final characterisation of the resulting critical points
\begin{theorem} {\bf (A)} The critical points of the $({\rm SU}(3),{\rm U}(2))$ MDM functional (\ref{lagr-gw-intr}) are four-dimensional constant scalar curvature almost-K\"ahler manifolds. 
\end{theorem}
Note that it is {\it not} implied that every constant scalar curvature almost-K\"ahler four-dimensional manifold can arise as a critical point of the MDM functional. 

We also remark that the functional (\ref{lagr-gw-intr}) could be considered in dimensions higher than four. The Euler-Lagrange equations (\ref{feqs-intr}) would be unchanged, and (\ref{main-feq-intr}) would be similar, with only differences in some numerical coefficients. What is very different in higher dimensions, however, is that $J$-invariance of Ricci no longer implies $d\rho=0$ for an almost-K\"ahler metric. The argument \cite{Draghici} that leads to this conclusion is crucially four-dimensional, see Section 3.7 where we spell it out. So, no constant scalar curvature can be deduced from the Euler-Lagrange equations for (\ref{lagr-gw-intr}) in higher dimensions. 

Another remark is that the specific relative coefficients in the Lagrangian in (\ref{lagr-gw-intr}) result from the MDM construction, or alternatively from the condition that the $\mathbb{CP}^2$ geometry with $\omega = \im da$ and $s=12$ is the critical point. If one removes this condition, then one could consider a more general family of functionals of the type (\ref{lagr-intr}) but with arbitrary $\lambda,\mu,\nu$. This is still a family of ${\rm U}(2)$ invariant functionals for an almost-Hermitian structure in four dimensions, as well as a ${\rm U}(2)$ connection. The conclusion that the Euler-Lagrange equations that follow for the ${\rm SU}(2)$ part $A$ of the connection $w$ are algebraic is still unchanged. Evaluating the functional on the solution $A={\bf w}$ still results in a functional of the type (\ref{lagr-gw-intr}) but this time with unrelated coefficients in front of the volume and $\omega \wedge \im da$ terms. This result in essentially the same Euler-Lagrange equations as (\ref{feqs-intr}), (\ref{main-feq-intr}), but with more general coefficients on the right-hand side of (\ref{main-feq-intr}). The conclusion that the critical points are constant scalar curvature almost-K\"ahler 4-manifolds is unchanged. 

The bottom line is that one can change the starting point, and instead of considering the MDM construction, which is only applicable in the context of $H$-structures of Cartan type $(G,H)$, consider the most general $H$-invariant functionals that can be built for the $H$ soldering form and an $H$-connection. In the case $H={\rm U}(2)$ considered here this still results in the functional (\ref{lagr-intr}), but with no condition on the parameters $\lambda,\mu,\nu$. This still results in the same, i.e. constant scalar curvature almost-K\"ahler 4-manifolds, critical points. 

All statements above are local. We can obtain further characterisations of the critical points by assuming that $M$ is a compact 4-manifold. The strongest result we can get uses the following theorem \cite{Draghici99}
\begin{theorem} {\bf (Draghici 99)} Let $(M,g,\omega)$ be a compact four-dimensional almost-K\"ahler manifold with J-invariant Ricci curvature and non-negative scalar curvature. Then the almost complex structure is integrable and $(M,g,\omega)$ is K\"ahler. 
\end{theorem}
Applying this result to our setup directly gives
\begin{corollary} {\bf (A)} Assume $M$ is compact 4-dimensional and the scalar curvature is non-negative. Then the critical points of the $({\rm SU}(3),{\rm U}(2))$ MDM functional (\ref{lagr-gw-intr}) are $s\geq 0$ constant scalar curvature K\"ahler (cscK) 4-manifolds.
\end{corollary}
Making an additional assumption about the first Chern class $c_1(TM,J)$ we get the following result
\begin{corollary} {\bf (B)} Assume $M$ is compact 4-dimensional, the scalar curvature is non-negative and  the first Chern class satisfies $2\pi c_1(TM,J)= \lambda [\omega]$ for some $\lambda\in \R$. Then the critical points of the MDM functional (\ref{lagr-gw-intr}) are $s\geq 0$ K\"ahler-Einstein 4-manifolds.
\end{corollary}
The argument here is standard in K\"ahler geometry. Assuming that the scalar curvature of a compact K\"ahler manifold is constant and the Ricci form is such that an Einstein metric can exist $[\rho]=\lambda [\omega]$ implies that the metric is K\"ahler-Einstein. This is shown by observing that $\rho$ in this case is harmonic, and being in the same cohomology class as $\omega$ must be a multiple of it. 

We can apply the same logic to our setup and prove the following generalisation which removes the $s\geq 0$ assumption
\begin{theorem} {\bf (B)} Assume $M$ is compact 4-dimensional and  the first Chern class satisfies $2\pi c_1(TM,J)= \lambda [\omega]$ for some $\lambda\in \R$.  The the critical points of the MDM functional (\ref{lagr-gw-intr}) are Einstein almost-K\"ahler 4-manifolds.
\end{theorem}
Indeed, in our almost-K\"ahler setup the logic is the same as in the more standard K\"ahler case. When the first Chern class $2\pi c_1(TM,J)$, which equals $[\im da]$, is a multiple of $[\omega]$, the equation (\ref{main-feq-intr}) implies that $[\rho]\sim[\omega]$. However, when the scalar curvature is constant and Ricci is J-invariant $\rho$ is a closed and a co-closed form, and thus necessarily $\rho=\lambda\omega$. 

Finally, we remind the reader that the Goldberg conjecture \cite{Goldberg} states that a compact Einstein almost-K\"ahler manifold is necessarily K\"ahler. This conjecture has been proven in the non-negative scalar curvature case by Sekigawa \cite{Sekigawa}, and \cite{Draghici99} strengthens this result by replacing the Einstein condition with a weaker $J$-invariance of Ricci condition. If Goldberg's conjecture holds then the only critical points of our functional in the setting $2\pi c_1(TM,J)= \lambda [\omega]$ would be K\"ahler-Einstein 4-manifolds.

The above statements demonstrate that the $({\rm SU}(3),{\rm U}(2))$ MDM construction produces an interesting and seemingly rich variational problem for almost-Hermitian structures in four dimensions. We expect that the same continues to be true for other possible setups. Of particular interest is the case $({\rm G}_2,{\rm SU}(3))$ in six dimensions, where the MDM construction gives an interesting functional for ${\rm SU}(3)$ structures. We hope to treat this case in a separate publication.

\section{SU(3) MDM functional}

\subsection{$\mathfrak{su}(2)$ description of the $\mathfrak{su}(3)$ connection}

We parametrise an $\mathfrak{su}(3)$-valued 1-form $\calA$ by its $\mathfrak{su}(2), \mathfrak{u}(1)$ and $\C^2$ components
\be 
\calA = \left( \begin{array}{cc} A + \id a & \Psi \\ - \Psi^\dagger & - 2a \end{array}\right).
\ee
Here $\calA$ is a $3\times 3$ anti-Hermitian matrix, $A$ is an $2\times 2$ anti-Hermitian matrix, $a\in \im \R$, and $\Psi\in \C^2$. All these fields are also 1-forms. The curvature of the SU(3) connection $\calA$ is given by
\be
\mathcal{F} = d\calA + \calA \wedge \calA,
\ee
which computes to be
\be
\mathcal{F} =\left( \begin{array}{cc} F + \id da - \Psi\wedge \Psi^\dagger & D\Psi + 3a\wedge \Psi \\ - (D\Psi+ 3a\wedge \Psi)^\dagger & - 2da - \Psi^\dagger\wedge \Psi\end{array}\right),
\ee
where
\be
F = dA + A\wedge A, \qquad D\Psi = d\Psi + A\wedge \Psi.
\ee
We also note that
\be
(D\Psi)^\dagger = D\Psi^\dagger = d\Psi^\dagger + \Psi^\dagger \wedge A.
\ee
In what follows we will often omit the wedge product sign. Thus, when not indicated, all products of differential forms are wedge products. 

\subsection{Pseudo-Pontryagin density}

The object 
\be
\int_M {\rm Tr}( \mathcal{F}\wedge \mathcal{F}) 
\ee
is a multiple of the topological invariant that only depends on the associated vector bundle (with fibres $\C^3$) on which $\calA$ is the connection. Motivated by what happens in MacDowell-Mansouri description of 4D GR, we deform this Lagrangian density by inserting under the trace a matrix that breaks ${\rm SU}(3)\to {\rm SU}(2)$. The natural such matrix is
\be
\gamma =  \left( \begin{array}{cc} \id & 0 \\ 0 & c \end{array}\right), \qquad c\in \R. 
\ee
Thus, we start by considering the theory 
\be\label{su3}
S[\calA] = \int_M {\rm Tr}( \gamma \mathcal{F} \wedge \mathcal{F}) .
\ee
We will later generalise this example to a functional of the type (\ref{functional-intr}) with more than one insertion. The functional computes to give
\be\nonumber
{\rm Tr}( \gamma \mathcal{F} \wedge \mathcal{F}) = {\rm Tr}\left( ( F+ \id da - \Psi\wedge \Psi^\dagger) \wedge ( F+ \id da - \Psi\wedge \Psi^\dagger)\right) \\ \nonumber
- (1+c) (D\Psi + 3a\wedge \Psi)^\dagger \wedge (D\Psi + 3a\wedge \Psi) + c (2da + \Psi^\dagger \wedge \Psi) \wedge (2da + \Psi^\dagger \wedge \Psi) .
\ee
We can now use the following integration by parts identities
\be
\frac{1}{2} D( \Psi^\dagger \wedge D\Psi + D\Psi^\dagger \wedge \Psi) = D\Psi^\dagger \wedge D\Psi - \Psi^\dagger F \Psi, \\ \nonumber
d( 3a \wedge \Psi^\dagger \wedge \Psi) = 3da \wedge \Psi^\dagger \wedge \Psi - 3a \wedge (D\Psi^\dagger \wedge \Psi - \Psi^\dagger \wedge D\Psi).
\ee
Here
\be
D\Psi^\dagger = d\Psi + \Psi^\dagger \wedge A = (D\Psi)^\dagger.
\ee
This means that the first term in the second line can be written as
\be
D\Psi^\dagger \wedge D\Psi + 3 (a\wedge \Psi)^\dagger \wedge D\Psi + 3 D\Psi^\dagger \wedge a \wedge \Psi + (a\wedge \Psi)^\dagger \wedge (a\wedge \Psi) = \\ \nonumber
D\Psi^\dagger \wedge D\Psi+ 3 a \wedge ( D\Psi^\dagger  \wedge \Psi - \Psi^\dagger \wedge D\Psi) = \\ \nonumber
 d( \frac{1}{2}\Psi^\dagger \wedge D\Psi + \frac{1}{2}D\Psi^\dagger \wedge \Psi -3a \wedge \Psi^\dagger \wedge \Psi)+ \Psi^\dagger F  \wedge \Psi   + 3 da \wedge \Psi^\dagger \wedge \Psi.
 \ee
 
 Using these identities, the action \eqref{su3} computes to 
 \begin{align}\nonumber
 S[\calA] &= \int_M {\rm Tr}(F\wedge F) + (2+4c) da\wedge da - (1+c) d( \frac{1}{2}\Psi^\dagger \wedge D\Psi + \frac{1}{2}D\Psi^\dagger \wedge \Psi -3a \wedge \Psi^\dagger \wedge \Psi)
 \\ \label{action-su3}
&+ (c-1) \left( {\rm Tr}( F \wedge \Psi\wedge \Psi^\dagger) + da\wedge \Psi^\dagger\wedge \Psi +  \Psi^\dagger \wedge \Psi \wedge  \Psi^\dagger \wedge \Psi \right),
 \end{align}
 where we have used
 \be
 {\rm Tr}(  \Psi\wedge \Psi^\dagger) = - \Psi^\dagger \wedge \Psi , \\ \nonumber
 {\rm Tr}( \Psi\wedge \Psi^\dagger \wedge  \Psi\wedge \Psi^\dagger) = - \Psi^\dagger \wedge \Psi \wedge  \Psi^\dagger \wedge \Psi.
 \ee
When $c=1$ the second line disappears, and the Lagrangian is a total derivative, as it should be for the Pontryagin density. When $c\not =1$, the first line is a total derivative, and the second line gives an interesting model.

\subsection{The arising model}

We now consider the theory with the action given by the second line in \eqref{action-su3} in its own right. It can be written as
\be\label{action-su3-model}
S[A,a,\Psi] = \int {\rm Tr}( (F - \id da) \wedge \Psi\wedge \Psi^\dagger) +  \Psi^\dagger \wedge \Psi \wedge  \Psi^\dagger \wedge \Psi.
\ee
The action is now viewed as a function of ${\rm SU}(2), {\rm U}(1)$ gauge fields $A,a$ and a $\C^2$-valued 1-form field $\Psi$. We note that $\Psi\wedge \Psi^\dagger$ is an anti-Hermitian $2\times 2$ matrix and so the action is real. 

The Euler-Lagrange equations are
\be\label{feqs}
D (\Psi\wedge \Psi^\dagger) \equiv d (\Psi\wedge \Psi^\dagger) + A\wedge (\Psi\wedge \Psi^\dagger) - (\Psi\wedge \Psi^\dagger)\wedge A = 0, \\ \nonumber
(F- \id da - 2  \Psi\wedge \Psi^\dagger) \wedge \Psi = 0. 
\ee
Here the tracefree and the trace parts of the first equation are those arising from the variation with respect to $A,a$ respectively, and the second line is the equation following upon the variation with respect to $\Psi^\dagger$. We note that the commutator in the first line is automatically tracefree, which means that the first equation contains the statement that $d (\Psi^\dagger\wedge \Psi)=0$. 

\subsection{Field equations directly}

It is instructive to also compute the field equations following from (\ref{su3}) directly, namely
\be
d_\mathcal{A} ( \gamma \mathcal{F}+ \mathcal{F}\gamma) =0.
\ee
We have
\be
\gamma \mathcal{F}+ \mathcal{F}\gamma = \left( \begin{array}{cc} 2(F+ \mathbb{I} da - \Psi\Psi^\dagger) & (1+c) (D\Psi+3a \Psi) \\ - (1+c)(D\Psi+ 3a\Psi)^\dagger & - 2c(2da + \Psi^\dagger \Psi) \end{array}\right).
\ee
A computation gives
\be \nonumber
&0= d_\mathcal{A} ( \gamma \mathcal{F}+ \mathcal{F}\gamma) =d (\gamma \mathcal{F}+ \mathcal{F}\gamma) + \mathcal{A}\wedge (\gamma \mathcal{F}+ \mathcal{F}\gamma)- (\gamma \mathcal{F} + \mathcal{F}\gamma) \wedge \mathcal{A} = \\
&(c-1) \left( \begin{array}{cc} D(\Psi \Psi^\dagger)  & ( F- \mathbb{I} da - 2 \Psi \Psi^\dagger) \Psi \\ \Psi^\dagger ( F- \mathbb{I} da - 2 \Psi \Psi^\dagger) & - d(\Psi^\dagger \Psi) \end{array}\right).
\ee
This vanishes when $c=1$ as it should. Note that the arising matrix is no longer tracefree, after the insertion of $\gamma$ which is not tracefree. The equations stored here are precisely the field equations (\ref{feqs}). 

\subsection{More general functionals}

The functional (\ref{action-su3-model}) has a very concrete set of coefficients in front of 
\be\label{blocks-su3}
\Psi^\dagger F \Psi, \quad da \Psi^\dagger \Psi, \quad (\Psi^\dagger\Psi)^2
\ee
terms. However, a simple generalisation of our construction (\ref{su3}) allows to produce a functional where these coefficients are more general. Indeed, consider the projectors on the $2\times 2$ and $1\times 1$ diagonal blocks
\be
P = \left( \begin{array}{cc} \mathbb{I} & 0 \\ 0 & 0 \end{array}\right), \qquad Q= \left( \begin{array}{cc} 0 & 0 \\ 0 & 1 \end{array}\right).
\ee
We can then consider the functional
\be\label{su3-general}
S[\mathcal{A}] = \int_M {\rm Tr}( (p P + q Q)\mathcal{F} (r P + s Q)\mathcal{F}).
\ee
This is precisely of the type we described in the Introduction (\ref{functional-intr}), with two different matrices $\Gamma_{1,2}$ inserted under the trace. 
It is clear that this functional is ${\rm U}(2)$-invariant. It is also clear that the whole construction is invariant under $(p,q)\to t (p,q), (r,s)\to t^{-1}(r,s)$, which means that we have effectively introduced only 3 parameters. Moreover, one of them can be removed by the overall rescaling of the action. So, we have introduced a 2-parameter family of functionals. It is clear that these parameters will translate into the relative coefficients between the 3 terms in (\ref{action-su3-model}). A simple computation gives
\be
S[\mathcal{A}]= \int_M {\rm Tr}\Big( pr ( F+ \id da - \Psi\wedge \Psi^\dagger)^2   \\ \nonumber
- (ps+qr) (D\Psi + 3a\wedge \Psi)^\dagger \wedge (D\Psi + 3a\wedge \Psi) + qs (2da + \Psi^\dagger \wedge \Psi)^2 \Big) .
\ee
Expanding and keeping only the non-topological terms, we get
\be\nonumber
S[\mathcal{A}]= \int_M {\rm Tr}\Big( pr (2\Psi^\dagger F\Psi + 2da \Psi^\dagger \Psi- (\Psi^\dagger \Psi)^2) - (ps+qr) (D\Psi^\dagger D\Psi + 3aD(\Psi^\dagger\Psi)) 
\\ \nonumber
+ qs (\Psi^\dagger \Psi)^2 + 4qs da \Psi^\dagger \wedge \Psi\Big).
\ee
Finally, integrating by parts we get
\be\label{su3-lagr-mu-nu}
S[\mathcal{A}]= \int_M {\rm Tr}\Big( \lambda  \Psi^\dagger F\Psi +  \mu da (\Psi^\dagger \Psi) + \nu  (\Psi^\dagger \Psi)^2\Big),
\ee
with
\be
\lambda = 2pr-ps-qr, \qquad \mu = 2pr+4qs-3(ps+qr), \qquad \nu = qs-pr.
\ee
This vanishes when $p,q,r,s=1$ as it should. Note, however, that the coefficients here are not independent, the following relation holds
\be\label{linear-relation}
3\lambda - \mu + 4\nu=0.
\ee
This proves our Proposition 1 from the Introduction. 

This above calculation shows that the freedom of choosing an ${\rm U}(2)$-invariant Lagrangian constructed from $\mathcal{F}$ is {\bf not} the same freedom as starting from the most general Lagrangian constructed as a linear combination of the ${\rm U}(2)$-invariant blocks (\ref{blocks-su3}). Indeed, the most general ${\rm U}(2)$-invariant Lagrangian would be of the form (\ref{su3-lagr-mu-nu}) with general coefficients, while we see that there is necessarily a linear relation between them. The same conclusion can be reached by analysing the field equations that follow (\ref{su3-lagr-mu-nu}). We will come back to this after we given the almost-Hermitian interpretation of the model. 

\section{Almost-Hermitian interpretation}

The purpose of this section is to explain why the functional (\ref{su3-lagr-mu-nu}) provides an action principle for almost-Hermitian geometry in 4D. 

\subsection{Metric and almost complex structure}

The field $\Psi$ is a $\C^2$-valued one form. We can write it more explicitly by parametrising
\be
\Psi = \left( \begin{array}{c} \alpha \\ \beta \end{array}\right).
\ee
Then 
\be
\Psi^\dagger \wedge \Psi = \bar{\alpha} \wedge \alpha  + \bar{\beta} \wedge \beta, \quad 
\Psi\wedge \Psi^\dagger = \left( \begin{array}{cc} \alpha \wedge \bar{\alpha} & \alpha \wedge \bar{\beta} \\ \beta \wedge \bar{\alpha} & \beta \wedge \bar{\beta}\end{array}\right).
\ee

We now extract the tracefree part of the matrix $\Psi\wedge \Psi^\dagger$, which is $\mathfrak{su}(2)$ Lie algebra valued
\be\label{Sigma}
\Sigma:= \Psi\wedge \Psi^\dagger\Big|_{\mathfrak{su}(2)} = \Psi\wedge \Psi^\dagger + \frac{\id}{2} \Psi^\dagger \wedge \Psi = \\ \nonumber
 \left( \begin{array}{cc}\frac{1}{2}( \alpha \wedge \bar{\alpha} - \beta\wedge \bar{\beta}) & \alpha \wedge \bar{\beta} \\ \beta \wedge \bar{\alpha} & -\frac{1}{2}( \alpha \wedge \bar{\alpha} - \beta\wedge \bar{\beta})\end{array}\right).
\ee
Let us also define a convenient multiple of the trace part
\be\label{omega}
\omega: = \frac{\im}{2} ( \alpha\wedge \bar{\alpha} + \beta\wedge \bar{\beta}) =  \frac{1}{2\im} \Psi^\dagger \Psi. 
\ee

There are now several equivalent ways of seeing both the metric and an orthogonal almost complex structure arising. The most direct way is to declare that the $\C$-valued 1-forms $\alpha,\beta$ are (1,0) forms for an almost complex structure that is compatible with the metric given by
\be\label{metric}
g = \alpha\odot \bar{\alpha} + \beta \odot \bar{\beta},
\ee
where $\odot$ stands for the symmetrised tensor product. The real 2-form $\omega$ is then the K\"ahler form. 

An alternative way to see the metric arising is to note that the triple of real 2-forms
\be
\frac{1}{2\im}( \alpha \wedge \bar{\alpha} - \beta\wedge \bar{\beta}) , \qquad {\rm Re}(\alpha \wedge \bar{\beta}), \quad {\rm Im}(\alpha \wedge \bar{\beta})
\ee
define a Riemannian signature metric, by requiring that this triple of 2-forms becomes anti-self-dual, and also requiring that the volume form is a multiple of $\alpha \wedge \bar{\alpha} \wedge \beta\wedge \bar{\beta}$. Throughout this article we will work in the orientation in which the K\"ahler 2-form becomes self-dual, so that the components of the matrix-valued 2-form $\Sigma$ are anti-self-dual. Note that $\Sigma\wedge \omega=0$. The metric defined by declaring the components of $\Sigma$ to span ASD 2-forms is (\ref{metric}). 

After the 2-form $\omega$ is declared to be the K\"ahler form, it directly defines an almost complex structure by raising one of its indices with the metric. It is not difficult to see that the resulting almost complex structure is the one whose action on the basis 1-forms is
\be
J (\alpha) = \im \alpha, \quad J(\beta) =\im \beta,
\ee
so that $g(\cdot,J\cdot) = \omega$. This is the same ACS as the one that arises by declaring $\Lambda^{(1,0)} = {\rm Span}(\alpha,\beta)$. 

To summarise, we have seen that the field $\Psi$ of our model defines almost-Hermitian geometry, in the sense of defining both a Reimannian signature metric (\ref{metric}) as well as a compatible orthogonal almost complex structure with the K\"ahler 2-form (\ref{omega}). We can also confirm this conclusion by a dimension count. The $\C^2$-valued 1-form $\Psi$ carries $4\times 4$ real functions. These are subject to ${\rm SU}(2)\times{\rm U}(1)$ gauge transformations. The ${\rm SU}(2)$ acts by mixing $\alpha,\beta$, while ${\rm U}(1)$ acts by
\be\label{u1-action}
(\alpha,\beta)\to (e^{\im\phi}\alpha,e^{-\im\phi}\beta). 
\ee
So, 4 out of 16 functions in $\Psi$ are gauge, leaving 12 functions. This is exactly the amount of information needed to specify a Riemannian metric (with its 10 components) and an almost complex structure (2 functions). We remind the reader that the space of almost-complex structures on $\R^4$, of a given orientation, can be described as the space of projective Weyl spinors of a given chirality. So, the space of almost-complex structures at a point is a copy of $\mathbb{CP}^1\sim S^2$. 

We can also arrive at the same conclusion in a different way. The group of orthogonal orientation-preserving transformations is ${\rm SO}(4)={\rm SU}(2)\times{\rm SU}(2)/\mathbb{Z}_2$. Let $\Lambda^\pm$ be the subspaces of self- and anti-self dual 2-forms (for a given metric). The two copies of ${\rm SU}(2)$ each acts on its own copy of $\Lambda^\pm$, while leaving the other space intact. Thus, let ${\rm SU}_-(2)$ be the group that acts non-trivially on the space $\Lambda^-$, while leaving all 2-forms in $\Lambda^+$ invariant. This is precisely the gauge group ${\rm SU}(2)$ of our model, if we identify 
\be
{\rm Span}( \Psi\wedge \Psi^\dagger \Big|_{\mathfrak{su}(2)} ) = \Lambda^-.
\ee
The trace part of the 2-form $\Psi\wedge \Psi^\dagger$, which is what we called $\omega$, see (\ref{omega}), is then in $\Lambda^+$.
The group ${\rm SU}_+(2)$ is the one that leaves invariant all 2-forms in $\Lambda^-$, and acts non-trivially on $\Lambda^+$. The subgroup of ${\rm SU}_+(2)$ that leaves $\omega$ invariant is the ${\rm U}(1)$ that acts as (\ref{u1-action}). This is the other factor in the gauge group of our model. What we have described is the standard ${\rm U}(2)={\rm SU}(2)\times{\rm U}(1)$ structure group of almost-Hermitian 4D geometry.

\subsection{More general field equations}

Now that we have the almost-Hermitian interpretation of the structure encoded by $\alpha,\beta$ at hand, we can analyse the implications of the field equations. 
It is clear that the most general functional (\ref{su3-general}) leads to the following field equations
\be\label{feqs-general}
D (\Psi\wedge \Psi^\dagger) \equiv d (\Psi\wedge \Psi^\dagger) + A\wedge (\Psi\wedge \Psi^\dagger) - (\Psi\wedge \Psi^\dagger)\wedge A = 0, \\ \nonumber
(\lambda F +\mu  \id da +2 \nu \Psi\wedge \Psi^\dagger) \wedge \Psi = 0,
\ee
where now $\lambda,\mu,\nu\in \R$ are parameters. From general considerations we know that these field equations must be satisfied by the model space $G/H$, which is when $\mathcal{F}=0$. It is instructive to verify this directly. 

Let us first consider what the vanishing of the off-diagonal part of the curvature $\mathcal{F}$ implies. This equation says that
\be\label{off-diagonal-F}
D\Psi + 3a\Psi=0.
\ee
Let us see what this implies for the exterior covariant derivative of $\Sigma$. We have
\be
D( \Psi \Psi^\dagger + \frac{\id}{2} \Psi^\dagger \Psi) = D\Psi \Psi^\dagger- \Psi D\Psi^\dagger + \frac{\id}{2} ( D\Psi^\dagger \Psi - \Psi^\dagger D\Psi) =
\\ \nonumber
- 3a\Psi \Psi^\dagger - \Psi (3a \Psi^\dagger) +  \frac{\id}{2}( 3a \Psi^\dagger \Psi - \Psi^\dagger (-3a \Psi)) = 0.
\ee
The equation $D\Sigma=0$ can be interpreted as an algebraic equation for the connection $A$. However, in the Appendix A we compute the exterior covariant derivative of the 2-form $\Sigma$ with respect to the ASD part ${\bf w}$ of the Levi-Civita connection, and we find $D_{\bf w}\Sigma=0$, see (\ref{d-Sigma-LC}). This means that the first equation in (\ref{feqs-general}) implies that $A={\bf w}$, where $\bf w$ is the chiral part of the Levi-Civita connection arising by decomposing $\mathfrak{so}(4)=\mathfrak{su}(2)\oplus \mathfrak{su}(2)$. Comparing (\ref{off-diagonal-F}) with (\ref{u2-torsion}), and taking into account $A={\bf w}$, we see that (\ref{off-diagonal-F}) implies that $3a=h$ and $T=0$. In particular, this means that the equation (\ref{off-diagonal-F}) implies that the torsion of the canonical ${\rm U}(2)$ connection is zero, which means that we are in the realm of K\"ahler geometry.

Overall, we see that (\ref{off-diagonal-F}) implies the first equation in (\ref{feqs-general}). It remains to see how the second equation in (\ref{feqs-general}) can be satisfied. We now take the diagonal parts of the $\mathcal{F}=0$ condition, which are
\be\label{diagonal-su3}
F + \mathbb{I} da = \Psi \Psi^\dagger, \qquad 2da + \Psi^\dagger \Psi =0.
\ee
Multiplying both of them by $\Psi$ (from the left) we get
\be
F \Psi + da \Psi = \Psi (\Psi^\dagger \Psi), \qquad da \Psi = -\frac{1}{2} \Psi (\Psi^\dagger \Psi) \quad \Rightarrow F\Psi = \frac{3}{2} \Psi (\Psi^\dagger \Psi).
\ee
Then the second field equation in (\ref{feqs-general}) just reduces to the linear relation (\ref{linear-relation}). We thus see that the relation (\ref{linear-relation}) must be satisfied in order for the field equations to hold on the configuration $\mathcal{F}=0$. The logic of the above derivation suggests that instead of considering the most general Lagrangian of the form (\ref{su3-general}) we can start with the most general ${\rm U}(2)$ invariant Lagrangian (\ref{su3-lagr-mu-nu}) provided we impose the condition on the coefficients that guarantees that $\mathcal{F}=0$ satisfies all the field equations. 

\subsection{Interpretation of the Cartan-flat geometries}

It can be expected that configurations that satisfy the Cartan-flat condition $\mathcal{F}=0$ correspond to the locally symmetric space $\mathbb{CP}^2={\rm SU}(3)/{\rm U}(2)$. Let us see this from the field equations. 

Let us note that, solving for $da$ in terms of $\Psi^\dagger\Psi$, the diagonal equations (\ref{diagonal-su3}) can also be rewritten as 
\be\label{diagonal-rewritten}
\im da= \omega, \qquad F=\Sigma,
\ee
where $\omega, \Sigma$ are defined by (\ref{omega}), (\ref{Sigma}) respectively. We now recall that the off-diagonal equations (\ref{off-diagonal-F}) allowed to identify $3a=h$ and $A={\bf w}$, where $h,{\bf w}$ are the parts of the Levi-Civita connection. This means that the ${\rm U}(1)$ connection $a$ gets identified with the Chern connection on the 
canonical line bundle. The first equation in (\ref{diagonal-rewritten}) then means that the curvature of this connection, which is the Ricci form, is a multiple of the K\"ahler form. This is the K\"ahler-Einstein condition. As before, we fix the orientation in which the K\"ahler form is self-dual $\omega\in \Lambda^+$. The second equation in (\ref{diagonal-rewritten}) is then a statement about the curvature of the ASD part of the Levi-Civita connection. It states that it is a constant multiple of the ASD 2-forms $\Sigma$, which also implies the Einstein condition, and in particular constant scalar curvature. But the equation (\ref{diagonal-rewritten}) also implies that the ASD part of the Weyl curvature is zero. This means that altogether the equations $\mathcal{F}=0$ imply that the geometry is K\"ahler-Einstein, with zero ASD part $W^-=0$ of the Weyl curvature tensor. In other words, we have a K\"ahler-Einstein self-dual 4-manifold. Such manifolds are known to be locally symmetric, see e.g. Lemma 7 of \cite{Derdzinski}.
Also the scalar curvature is strictly positive. This means that the only solution is the locally symmetric complex surface $\mathbb{CP}^2$. This is the expected solution, because we are considering a geometry with flat Cartan connection. This is the geometry with the total space of the ${\rm SU}(2)$ bundle over it identifiable with ${\rm SU}(3)$, which means $\mathbb{CP}^2$. 

\subsection{Functional for almost-Hermitian structures}

It is becoming clear that the Lagrangian (\ref{su3-lagr-mu-nu}) with (\ref{linear-relation}) satisfied describes almost-Hermitian geometry. When the fields are taken to satisfy $\mathcal{F}=0$ the torsion of this almost-Hermitian geometry (defined and discussed in the Appendix A) vanishes, and we are dealing with K\"ahler geometry. As we just discussed, the $\mathcal{F}=0$ condition also implies Einstein with positive scalar curvature and $W^-=0$, which is extremely restricting. Let us now understand what the Euler-Lagrange equations of the model imply more generally. 

Let us start by rewriting the Lagrangian (\ref{su3-lagr-mu-nu}) in metric terms. We have seen that the $A$ field equation implies $A={\bf w}$. Substituting this into the Lagrangian we get a functional depending solely on $\Psi$
\be
S[\Psi] = \int_M {\rm Tr}( \lambda \Psi^\dagger F({\bf w}) \Psi + \mu da (\Psi^\dagger \Psi) + \nu (\Psi^\dagger \Psi)^2).
\ee
In Appendix we compute the quantity $\Psi^\dagger F({\bf w}) \Psi$ and show that it equals (\ref{relation-pfp-R}) to minus the Ricci scalar times the volume form. The relation (\ref{omega-2}) gives the expression for $(\Psi^\dagger \Psi)^2$ in terms of the volume form. Thus, we can rewrite the arising functional as that of a compatible pair $(g,\omega)$.
Introducing the rescaled parameters $\tilde{\nu}=\nu/\lambda, \tilde{\mu}=\mu/\lambda$, and taking into account that by (\ref{linear-relation}) these parameters satisfy
\be\label{nu-mu}
4\tilde{\nu} = \tilde{\mu}-3,
\ee
the action takes the form
\be\label{act-g-w}
S[g,\omega] = -  \lambda \int_M dV ( s + 2\tilde{\mu}-6) -2\im \tilde{\mu} \omega \, da .
\ee
This is a functional of a Riemannian metric $g$ together with a compatible with it non-degenerate 2-form $\omega$, where the compatibility condition means that $J_\mu{}^\nu=\omega_{\mu\alpha} g^{\alpha\nu}$ is an almost complex structure, which can be written as
\be\label{g-w-compat}
\omega_{\mu\alpha} g^{\alpha\beta} \omega_{\beta\nu} = - g_{\mu\nu}.
\ee
Note that this is essentially the Einstein-Hilbert functional, but supplemented with the Lagrange multiplier that forces the K\"ahler form $\omega$ to be closed, as well as the condition (\ref{g-w-compat}) that constraints the variations of $g,\omega$. This proves Proposition 2 from the Introduction. 

\subsection{Consequences of the $\omega$ and $g$ compatibility}

If not for the relation (\ref{g-w-compat}) that must hold between the 2-form $\omega$ and the metric $g$, the Euler-Lagrange equations that would follow from (\ref{act-g-w}) would be a decoupled system of Einstein equations for the metric, as well as the condition that $\omega$ is closed. However, the fields $g,w$ are not free to vary independently, and this modifies the Euler-Lagrange equations.

Let us decipher the consequences of (\ref{g-w-compat}). The linearisation of this condition can be written as
\be\label{lin-compat}
J_\mu{}^\alpha \delta\omega_{\alpha\nu} + J_\nu{}^\alpha \delta\omega_{\alpha\mu} + J_\mu{}^\alpha J_\nu{}^\beta \delta g_{\alpha\beta} + \delta g_{\mu\nu}=0.
\ee
Here we have used the expression for the operator of the almost complex structure $\omega_{\mu\alpha} g^{\alpha\nu} = J_\mu{}^\nu$. To see what this condition implies we must decompose the variations $\delta \omega, \delta g$ into their ${\rm U}(2)$ irreducible parts. We have
\be
\delta \omega \in \Lambda_\R^{1,1} \oplus {\rm Re}(\Lambda^{2,0}).
\ee
A simple counting is useful here. There are six real components of $\delta \omega$. The real dimension of $\Lambda_\R^{1,1}$ is four. Indeed, in a basis $\alpha,\beta$ of $\Lambda^{1,0}$ the basis for $\Lambda_\R^{1,1}$ is $(1/\im)\alpha\bar{\alpha}, (1/\im)\beta\bar{\beta}$, as well as the real and imaginary parts of $\alpha\bar{\beta}$. The space $\Lambda^{2,0}$ has complex dimension one, so real dimension two. 

Now the first two terms in (\ref{lin-compat}) take the object $\delta \omega(J\cdot,\cdot)$ and symmetrise it in its two slots. Applying this operation to the part of $\delta\omega$ in $\Lambda^{2,0}$ produces zero, as $\delta \omega(J\cdot,\cdot)$ in this case is anti-symmetric. So, only the $\Lambda^{1,1}_\R$ part of $\delta\omega$ survives in (\ref{lin-compat}). 

Let us similarly decompose $\delta g$. We can write
\be
S^2 T^*M = S^{1,1}_\R \oplus {\rm Re}( S^{2,0}). 
\ee
The dimension count here is again useful. The dimension of $S^{1,1}_\R$ is four. The complex dimension of $S^{2,0}$ is three, which gives real dimension six. Together with four real components in $S^{1,1}_\R$ this gives the ten components of $\delta g$. 

Now, the two last terms in (\ref{lin-compat}) are of the schematic form $\delta g(J\cdot,J\cdot) + \delta g(\cdot,\cdot)$. It is clear that this projects out the $S^{2,0}$ part of $\delta g$, and only the $S^{1,1}_\R$ part survives. This shows that (\ref{lin-compat}) says
\be\label{lin-conds}
\delta \omega^{1,1}(J\cdot,\cdot) + \delta g^{1,1} =0.
\ee
In words, for the 2-form $\omega$ and a metric $g$ to satisfy (\ref{g-w-compat}) their $(1,1)$ parts must be related to each other in a variation, while the $\Lambda^{2,0}$ part of $\delta w$ and the $S^{2,0}$ part of $\delta g$ are unconstrained. The relation (\ref{lin-conds}) is four real conditions that are imposed on the 16-dimensional space $\delta g,\delta \omega$, leaving the variational space of dimension twelve, which is the correct dimension of the space of almost-Hermitian structures in four dimensions. 

\subsection{Euler-Lagrange equations}

We can now determine the Euler-Lagrange equations following by minimising (\ref{act-g-w})). 
The equation obtained by varying the action with respect to $a$ is simply
\be
d\omega =0.
\ee
Thus, generally, our model is about almost-K\"ahler geometry, see e.g. \cite{apostolov} for a survey. 

To obtain the EL equations that follow by minimising (\ref{act-g-w}) with respect to $\omega, g$ we recall that the variations with respect to $\delta w^{2,0}$ and $\delta g^{2,0}$ are independent. We also use the standard fact that the variational derivative of the Einstein-Hilbert functional (the first term in (\ref{act-g-w})) is $Ric- (s/2) g$. This immediately gives
\be
Ric^{2,0}=0, \qquad da^{2,0}=0.
\ee
Alternatively, we can say that the critical points of the functional (\ref{act-g-w}) are metrics with Ricci tensor that is J-invariant, and 1-forms $a$ that have J-invariant exterior derivative
\be
Ric(J\cdot,J\cdot)= Ric(\cdot,\cdot), \qquad da(J\cdot,J\cdot) = da(\cdot,\cdot).
\ee
Yet another equivalent way of stating these conditions is
\be
Ric\in S^{1,1}_\R, \qquad da\in \Lambda^{1,1}_\R.
\ee

The obtained condition on the Ricci tensor is exactly the same as follows from a related but different variational principle considered in \cite{Blair-Ianus}. In this reference, the authors study the Euler-Lagrange equations arising from the functional given by the integrated scalar curvature, but restricted to metrics compatible with a given symplectic 2-form. Because only the $(2,0)$ part of the metric is allowed to vary in this case, they find that the critical points are metrics whose Ricci tensor is $J$-invariant. Our functional (\ref{act-g-w}) is different from the one considered in \cite{Blair-Ianus} because we allow both $(g,\omega)$ to vary (while remaining compatible), and there is also an explicit Lagrange multiplier field in the Lagrangian imposing the condition that $d\omega=0$. 

The fact that the Lagrangian (\ref{act-g-w}) depends on both fields $(g,\omega)$ whose $(1,1)$ variations are related means that, unlike in \cite{Blair-Ianus}, we also get an interesting EL equation for the $(1,1)$ part of the Ricci tensor. It is clear that the arising EL equation simply says that the $(1,1)$ parts of Ricci and $da$ are proportional to each other, plus a contribution (proportional to $\omega$) and coming from the variation of the volume form. We get an equation that can be written in terms of 2-forms
\be\label{main-feq}
\rho = \frac{1}{2}\omega (s -6)  +\tilde{\mu} (\omega - \im da).
\ee
Here we have used (\ref{nu-mu}) and the fact that $Ric$ is $(1,1)$ and therefore $\rho=Ric(J\cdot,\cdot)$ is a 2-form in $\Lambda^{1,1}_\R$. We have also assumed $\lambda\not=0$ and divided by it. Note that all terms in this equations are in $\Lambda^{1,1}_\R$. The equation (\ref{main-feq}) is satisfied by $\mathbb{CP}^2$ which corresponds to $\omega=\im da$ and $s=12$, and for which we also have $\rho=s/4 \omega$. This finishes the proof of Proposition 3 from the Introduction. 

\subsection{Consequences of the field equations}

Taking the exterior derivative of (\ref{main-feq}) we have
\be
d\rho = \frac{1}{2} ds \wedge \omega.
\ee
On the other hand, \cite{Draghici} showed that any almost-K\"ahler 4-dimensional manifold with J-invariant Ricci curvature satisfies $d\rho=0$. The argument for this is very simple. The differential Bianchi identity for the Ricci tensor states that the symmetric rank two tensor $Ric - (s/2)g$ is co-closed. If the Ricci tensor is $J$-invariant $Ric(J\cdot, J\cdot)=Ric(\cdot,\cdot)$, then a simple computation, see a lemma in \cite{Draghici} shows that the corresponding 2-form $\rho- (s/2)\omega$ is co-closed as well. One then has
\be
0=d\star (\rho - \frac{s}{2}\omega) = d \star ( \rho - \frac{s}{4}\omega - \frac{s}{4}\omega) = d( - \rho + \frac{s}{4}\omega - \frac{s}{4}\omega) = - d\rho,
\ee
where we have used the fact that in four dimensions the self- and anti-self dual parts of the 2-form $\rho$ are given by
\be
\rho = \rho^+ + \rho^-, \qquad \rho^+ = \frac{s}{4} \omega, \qquad \rho^- = \rho -  \frac{s}{4} \omega.
\ee

This shows that $d\rho=0$ and thus $ds\wedge \omega=0$, which implies $ds=0$. Thus, as stated in Theorem A of the Introduction, the critical points of our functional are constant scalar curvature almost-K\"ahler metrics (cscAK). The arguments proving Theorems (B), (C) are straightforward, and are given already in the Introduction, so will not be repeated here. 

\section*{Acknowledgements} KK is grateful to Joel Fine for a discussion on the topics of this paper, and to Vestislav Apostolov for correspondence. PA acknowledges support from ANID-FONDECYT Regular grant No. 1230112.

\appendix

\section{Almost Hermitian geometry in four dimensions}

The purpose of this Appendix is to derive some standard facts about almost-complex geometry in four dimensions that we use in the main text. Our desire is to be as explicit as possible. We achieve this by doing all the calculations using a conveniently adapted frame.  

\subsection{Connection and intrinsic torsion of an almost-Hermitian structure}

Let $e^{1,2,3,4}$ be an orthonormal coframe in four dimensions. Given an orthogonal almost complex structure, let $\alpha,\beta$ space the space of $(1,0)$ forms. Is is natural to restrict to coframes adapted to the almost complex structure, so that
\be
\alpha = e^1 + \im e^2, \qquad \beta = e^3+\im e^4.
\ee
The fundamental theorem of Riemannian geometry is that there exists a unique metric and torsion free connection, which we call $w^I{}_J$. It is metric $w^{IJ}=w^{[IJ]}$, and the condition of zero torsion is the Cartan's first structure equation
\be
d e^I = - w^I{}_J e^J.
\ee
Let us compute $d\alpha, d\beta$. A computation gives
\be\nonumber
d\alpha = \im w^{12} \alpha + \frac{\beta}{2}( -(w^{24} +\im w^{23}) + \im (w^{14}+\im w^{13})) + \frac{\bar{\beta}}{2} ( (w^{24}-\im w^{23}) - \im (w^{14}-\im w^{13})),
\\ \nonumber
d\beta = \im w^{34} \beta + \frac{\alpha}{2}( (w^{24} -\im w^{23}) + \im (w^{14}-\im w^{13})) + \frac{\bar{\alpha}}{2} ( -(w^{24}-\im w^{23}) + \im (w^{14}-\im w^{13})).
\ee
Let us introduce
\be\label{w3-w}
w^3 = -w^{12}+w^{34}, \quad  w = -(w^{24} +\im w^{23}) + \im (w^{14}+\im w^{13}), \\ \nonumber 
\tilde{w}^3 = -w^{12}-w^{34}, \quad T= (w^{24}-\im w^{23}) - \im (w^{14}-\im w^{13}).
\ee
Then, introducing
\be
\Psi = \left( \begin{array}{c} \alpha \\ \beta \end{array}\right), \quad {\bf w} = \frac{1}{2} \left(\begin{array}{cc} \im w^3 & w \\ - \bar{w} & - \im w^3\end{array}\right), \quad h= \frac{\im}{2} \tilde{w}^3, \quad 
\quad \epsilon = \left( \begin{array}{cc} 0 & 1 \\ -1 & 0 \end{array}\right),
\ee
we can rewrite the formulas for $d\alpha,d\beta$ in matrix form
\be\label{u2-torsion}
d \Psi + ({\bf w}+  \mathbb{I} h) \Psi = - \frac{1}{2}T \epsilon \bar{\Psi}.
\ee
The terms of the left-hand side are the covariant exterior derivative of $\Psi$ with respect to the canonical ${\rm U}(2)$ connection, which is the ${\rm U}(2)$ part of the Levi-Civita connection. The terms on the right-hand side encode the so-called intrinsic torsion of the ${\rm U}(2)$ structure. If $T$ is non-zero, the ${\rm U}(2)$ structure is not parallel with respect to the Levi-Civita connection, and hence not integrable. On the other hand, if $T$ vanishes, the almost-Hermitian structure is integrable, which means that the geometry is K\"ahler. The equation (\ref{u2-torsion}) is standard in the literature, see e.g. equation (2.57) in \cite{Bryant-6D}.

\subsection{Useful consequences of the ${\rm U}(2)$ structure equation}

We will refer to (\ref{u2-torsion}) as the ${\rm U}(2)$ structure equation. Let us use this equation to obtain some useful consequences. We first apply (\ref{u2-torsion}) to compute 
\be
D_{\bf w} (\Psi\Psi^\dagger) = D_{{\bf w}+\mathbb{I} h} (\Psi\Psi^\dagger) = D_{{\bf w}+\mathbb{I} h} \Psi \Psi^\dagger - \Psi D_{{\bf w}+\mathbb{I} h} \Psi^\dagger = 
\\ \nonumber
- \frac{1}{2} T \epsilon \bar{\Psi} \Psi^\dagger + \frac{1}{2} \bar{T} \Psi \Psi^T \epsilon.
\ee
We have
\be
\epsilon \bar{\Psi} \Psi^\dagger = \left( \begin{array}{c} \bar{\beta} \\ -\bar{\alpha} \end{array}\right) \left( \begin{array}{cc} \bar{\alpha} & \bar{\beta} \end{array}\right)=
- \bar{\Omega} \mathbb{I}, \qquad \bar{\Omega} := \bar{\alpha}\bar{\beta},
\\ \nonumber
\Psi \Psi^T \epsilon = \left( \begin{array}{c} \alpha \\ \beta \end{array}\right) \left( \begin{array}{cc} -\beta & \alpha \end{array}\right)=
- \Omega \mathbb{I}, \qquad \Omega := \alpha\beta.
\ee
and so
\be
D_{\bf w} (\Psi\Psi^\dagger) =  \frac{1}{2} T \bar{\Omega} \mathbb{I} -\frac{1}{2} \bar{T} \Omega \mathbb{I} .
\ee
We also have for the trace part 
\be
d \Psi^\dagger \Psi = D_{{\bf w}+\mathbb{I} h} (\Psi^\dagger \Psi) = \frac{1}{2} \bar{T} \Psi^T \epsilon \Psi - \frac{1}{2} T \Psi^\dagger  \epsilon \bar{\Psi},
\ee
which we rewrite using
\be
\Psi^T \epsilon \Psi = 2\Omega
\ee
as
\be
d \Psi^\dagger \Psi =  \bar{T} \Omega - T \bar{\Omega}.
\ee
If we introduce
\be\label{d-omega-app}
\omega = \frac{1}{2} (\alpha \bar{\alpha} + \beta\bar{\beta}) = - \frac{1}{2} \Psi^\dagger \Psi ,
\ee
we have
\be\label{d-omega}
d\omega = - \frac{1}{2} \bar{T} \Omega +  \frac{1}{2} T \bar{\Omega}.
\ee

If we introduce the tracefree part of the matrix $\Psi\Psi^\dagger$
\be
\Sigma := \Psi\Psi^\dagger + \frac{1}{2} \mathbb{I} \Psi^\dagger \Psi, 
\ee
we have
\be\label{d-Sigma-LC}
D_{\bf w} \Sigma = \frac{1}{2} T \bar{\Omega} \mathbb{I} -\frac{1}{2} \bar{T} \Omega \mathbb{I} + \frac{1}{2} \mathbb{I} (\bar{T} \Omega - T\bar{\Omega})=0.
\ee
 Thus, the tracefree part of the matrix $\Psi\Psi^\dagger$ is covariantly constant with respect to the connection $\bf w$ (or ${\bf w}+\mathbb{I} h$), while the trace part 2-form $\omega$ satisfies the equation (\ref{d-omega}). 
 
 Let us also compute the exterior derivatives of $\Omega,\bar{\Omega}$. We have
 \be
 2d\Omega = d( \Psi^T \epsilon \Psi) = D_{\bf w} \Psi^T \epsilon \Psi - \Psi^T \epsilon D_{\bf w} \Psi. 
 \ee
 Here we can rewrite everything using the exterior covariant derivative with respect to $\bf w$, because for any tracefree matrix ${\bf w}$ we have ${\bf w}^T \epsilon = -\epsilon {\bf w}$, and
 so
 \be
 D_{\bf w} \Psi^T \epsilon \Psi - \Psi^T \epsilon D_{\bf w}\Psi = (d\Psi^T - \Psi^T {\bf w}^T) \epsilon \Psi - \Psi^T \epsilon (d\Psi + {\bf w}\Psi) = \\ \nonumber
 d\Psi^T  \epsilon \Psi - \Psi^T \epsilon d\Psi = d( \Psi^T \epsilon \Psi) .
 \ee
 Continuing with the calculation
 \be
 2d\Omega = (-\frac{1}{2} T \Psi^\dagger \epsilon^T - h \Psi^T) \epsilon \Psi - \Psi^T \epsilon ( -\frac{1}{2} T \epsilon \bar{\Psi} - h\Psi) =
 2T \omega - 2a \Omega,
 \ee
 or, equivalently
 \be\label{d-Omega}
 (d+h)\Omega = T \omega.
 \ee
Together the formulas (\ref{d-omega}) and (\ref{d-Omega}) show that when the intrinsic torsion $T$ vanishes the 2-form $\omega$ is closed and $(d+h)\Omega=0$. The second equation implies intregrability of the almost-complex structure, and the first then implies the K\"ahler condition. So, once again, we see that the vanishing of the intrinsic torsion $T$ implies K\"ahler. 

The calculations we performed are standard in the context of $H$-structures. One first establishes an analog of the formula (\ref{u2-torsion}) for the exterior covariant derivative of the adapted frame with respect to the canonical connection. The right-hand side of such a formula is the intrinsic torsion of the $H$-structure. One can then use this formula to compute the exterior derivatives of the ${\rm U}(2)$ invariant (or covariant, as in the case of $\Omega$) differential forms. These do not involve the canonical connection, and become expressions involving only the intrinsic torsion. Equations of the type (\ref{d-omega}) and (\ref{d-Omega}) could alternatively be derived without (\ref{u2-torsion}), using the decomposition of the covariant derivatives of $\omega,\Omega$ into irreducible representations of the group ${\rm U}(2)$. But the presented here direct derivation is simple and instructive. 

\subsection{Curvature of the ${\rm SU}(2)$ connection $\bf w$}

We now explicitly compute the curvature of the connection $\bf w$. We have
\be
F({\bf w}) = d{\bf w} + {\bf w}\wedge {\bf w} = \frac{1}{2} \left( \begin{array}{cc} \im F^3 & F \\ - \bar{F} & - \im F^3\end{array}\right),
\ee
where
\be
F^3 = dw^3 +\frac{\im}{2} w\wedge \bar{w}, \qquad F= dw + \im w^3 \wedge w.
\ee
A simple computation gives
\be
\Psi^\dagger F \Psi = - \frac{\im}{2} (\alpha\bar{\alpha} - \beta\bar{\beta}) F^3 + \frac{1}{2} \alpha\bar{\beta} \bar{F} + \frac{1}{2}  \bar{\alpha} \beta F.
\ee
Our task now is to compute the Ricci scalar (times the volume form $dV$), and compare this to $\Psi^\dagger F\Psi$. 

On one hand, recalling the expressions (\ref{w3-w}) for $w^3, w$ in terms of the components of the Levi-Civita connection, we have
\be\nonumber
F^3= d( - w^{12}+w^{34}) - w^{23} w^{24} - w^{13} w^{14} - w^{24} w^{14} - w^{24} w^{13}, \\ \nonumber
F = d( - (w^{24}+\im w^{23}) +\im (w^{14}+\im w^{13})) + \im ( - w^{12}+w^{34})( - (w^{24}+\im w^{23}) +\im (w^{14}+\im w^{13})).
\ee

On the other hand, we have
\[
s dV  = \frac{1}{4} \epsilon_{IJKL} R^{IJ} e^K e^L = R^{12} e^{34} - R^{13} e^{24} + R^{14} e^{23} + R^{23} e^{14} - R^{24} e^{13} + R^{34} e^{12},
\]
where $e^{IJ} = e^I e^J$ and
\be
R^{IJ} = dw^{IJ} + w^I{}_K w^{KJ}.
\ee
However, the expression 
\be
\frac{1}{2} R^{IJ} e_{IJ} = R^{12} e^{12} + R^{13} e^{12} + R^{14} e^{14} + R^{23} e^{23} + R^{24} e^{24} + R^{34} e^{34}
\ee
vanishes because of the (algebraic) Bianchi identity, which holds true when $w^{IJ}$ is the Levi-Civita connection. So, we can obtain an alternative expression for the Ricci scalar subtracting the zero
\[
s dV =   (R^{34}-R^{12}) (e^{12}-e^{34})- (R^{13} +R^{24})(e^{24}+e^{13}) + (R^{14} -R^{23})(e^{23} -e^{14}).
\]
Let us now compute the 2-form combinations that appear here. 
Using $\alpha = e^1+\im e^2, \beta = e^3+\im e^4$ we have
\[
e^{12} - e^{34}= \frac{\im}{2} (\alpha\bar{\alpha} - \beta\bar{\beta}), \quad e^{13} +e^{24} = \frac{1}{2}(\alpha\bar{\beta}+ \bar{\alpha}\beta), \quad e^{23} -e^{14} = \frac{1}{2\im}(\alpha\bar{\beta}- \bar{\alpha}\beta).
\]

We also have
\be\nonumber
R^{12} = dw^{12} - w^{13} w^{23} - w^{14} w^{24}, \quad R^{13} = dw^{13} + w^{12} w^{23} - w^{14} w^{34}, \\ \nonumber
R^{14} = dw^{14} + w^{12} w^{24} + w^{13} w^{34},\quad R^{23} = dw^{23} - w^{12} w^{13} - w^{24} w^{34}, \\ \nonumber 
R^{24} = dw^{24} - w^{12} w^{14} + w^{23} w^{34}, \quad R^{34} = dw^{34} - w^{13} w^{14} - w^{23} w^{24},
\ee
and therefore
\be\nonumber
R^{34}-R^{12} = d(w^{34}-w^{12}) + w^{13} w^{23} + w^{14} w^{24}- w^{13} w^{14} - w^{23} w^{24}= dw^3 + \frac{\im}{2} w\bar{w} = F^3.
\ee
The coefficients of $\alpha\bar{\beta},\bar{\alpha}\beta$ in $s dV$ are then
\[
- \frac{1}{2}(R^{13} +R^{24}) + \frac{1}{2\im}(R^{14} -R^{23})= \frac{1}{2} \bar{F}, \quad
- \frac{1}{2}(R^{13} +R^{24}) - \frac{1}{2\im}(R^{14} -R^{23})= \frac{1}{2} F.
\]
Thus, we have the following relation between the quantity $\Psi^\dagger F\Psi$ and the Ricci scalar
\be\label{relation-pfp-R}
\Psi^\dagger F\Psi = - s dV.
\ee
Another useful relation is
\be\label{omega-2}
(\Psi^\dagger \Psi)^2 = 2 \alpha\bar{\alpha}\beta \bar{\beta} = - 8 e^{1234} = - 8 dV.
\ee


\begin{thebibliography}{99}

\bibitem{MacDowell:1977jt}
S.~W.~MacDowell and F.~Mansouri,
``Unified Geometric Theory of Gravity and Supergravity,''
Phys. Rev. Lett. \textbf{38}, 739 (1977)
[erratum: Phys. Rev. Lett. \textbf{38}, 1376 (1977)]
doi:10.1103/PhysRevLett.38.739

\bibitem{Wise:2006sm}
D.~K.~Wise,
``MacDowell-Mansouri gravity and Cartan geometry,''
Class. Quant. Grav. \textbf{27}, 155010 (2010)
doi:10.1088/0264-9381/27/15/155010
[arXiv:gr-qc/0611154 [gr-qc]].

\bibitem{Chamseddine:1990gk}
A.~H.~Chamseddine,
``Topological gravity and supergravity in various dimensions,''
Nucl. Phys. B \textbf{346}, 213-234 (1990)
doi:10.1016/0550-3213(90)90245-9

\bibitem{Castellani:2017vbi}
L.~Castellani,
``A locally supersymmetric $SO(10,2)$ invariant action for $D=12$ supergravity,''
JHEP \textbf{06}, 061 (2017)
doi:10.1007/JHEP06(2017)061
[arXiv:1705.00638 [hep-th]].

\bibitem{Sharpe} R.\ W.\ Sharpe, "Differential geometry: Cartan's generalisation of Klein's Erlangen programme," Springer, 1997. 

\bibitem{HH} F.\ Hirzebruch and H.\ Hopf, "Felder von Fl\"achenelementen in 4-dimensionalen Mannigfaltigkeiten," Math. Ann., Vol. 136 156–172 (1958).

\bibitem{Draghici} T.\ Draghici, "On some 4-dimensional almost K\"ahler manifolds," Kodai Math. J. 18 (1995), 156–163.

\bibitem{Draghici99} T.\ Draghici, "Almost K\"ahler 4-manifolds with J-invariant Ricci tensor," Houston J. Math. Vol. 25 no. 1 133–145 (1999). 

\bibitem{Goldberg} S.\ I.\ Goldberg, "Integrability of almost-K\"ahler manifolds," Proc. Amer. Math. Soc. Vol. 21
96–100 (1969).

\bibitem{Sekigawa} K.\ Sekigawa, "On some compact Einstein almost-K\"ahler manifolds," J. Math. Soc. Japan Vol. 36
677–684 (1987). 

\bibitem{Derdzinski} A.\ Derdzinski, "Self-dual Kähler manifolds and Einstein manifolds of dimension four,"
Compositio Mathematica, tome 49, no 3 (1983), p. 405-433.

\bibitem{apostolov} V.\ Apostolov and T.\ Draghici, "The curvature and the integrability of almost-K\"ahler manifolds: A survey," in *Symplectic and contact topology: interactions and perspectives*, vol. 35, pp. 25-53, Amer. Math. Soc., 2003, arXiv:math/0302152 [math.DG].

\bibitem{Blair-Ianus} D.\ Blair and S.\ Ianus, "Critical asssociated metrics on symplectic manifolds," Contemporary Mathematics, Vol. 51 23-29 (1986).

\bibitem{Bryant-6D} R.\ Bryant, "ON THE GEOMETRY OF ALMOST COMPLEX 6-MANIFOLDS," math.DG/0508428.



\end{thebibliography}
\end{document}